\documentclass[openany,letter]{article}

\let\svthefootnote\thefootnote
\usepackage{multirow}
\usepackage{enumerate}
\usepackage[sumlimits]{amsmath}
\usepackage{amsfonts,amssymb}
\usepackage{enumitem}

\usepackage[margin=2.5cm]{geometry}
\usepackage[english]{babel}
\usepackage[all]{xy}
\usepackage{graphics}
\usepackage{mathrsfs}
\usepackage{ntheorem}
\usepackage{shuffle}

\usepackage{makeidx}

\usepackage[compact]{titlesec}
\usepackage{fancyhdr}

\xyoption{import}

\makeindex

\fancypagestyle{plain}{%
\fancyhf{} 
\fancyfoot[C]{\thepage} 

}
\titleformat{\chapter}[display]
{\bfseries\Huge} {\vspace{-2cm}\filleft\Huge \chaptertitlename
\hspace{.5cm}\thechapter} {2ex} {\titlerule[1pt] \vspace{1pt}
\titlerule
\vspace{1ex}%
\filright}
[\vspace{1ex}%
\titlerule]

\titleformat{\section}[hang]
{\Large\sffamily}
{\thesection.--}
{2pt}
{\vspace{2pt}\Large}

\newcommand{\fin}{\hfill$\blacksquare$}
\newcommand{\enf}[1]{\textsf{\textbf{#1}}}
\newcommand{\dob}[1]{\mathbb{#1}}

\numberwithin{equation}{section}

\theoremstyle{plain}
\theoremheaderfont{\normalfont\bfseries}
\theorembodyfont{\slshape}
\theoremseparator{.}
\newtheorem{teo}{Theorem}[section]
    \theoremclass{teo}
    \theoremstyle{break}
    
    \theoremclass{teo}
    
    \theoremclass{teo}
    \theoremstyle{break}
    
    \theoremclass{teo}
    \theoremstyle{nonumberplain}

\theoremstyle{plain}
\theoremheaderfont{\normalfont\bfseries}
\theorembodyfont{\slshape}
\theoremseparator{.}
\newtheorem{lema}[teo]{Lemma}
    \theoremclass{lema}

\theoremstyle{plain}
\theoremheaderfont{\normalfont\bfseries}
\theorembodyfont{\slshape}
\theoremseparator{.}
\newtheorem{p}[teo]{Proposition}
    \theoremclass{p}
    \theoremstyle{break}
    
    \theoremclass{p}
    
    \theoremclass{p}
    \theoremstyle{break}

\theoremstyle{plain}
\theoremheaderfont{\normalfont\bfseries}
\theorembodyfont{\slshape}
\theoremseparator{.}
\theoremstyle{change}
\theorembodyfont{\upshape}
\theoremseparator{:}
\newtheorem{ej}[teo]{Example}
    \theoremclass{ej}
    \theoremstyle{changebreak}

\theoremstyle{nonumberplain}
\theorembodyfont{\upshape}
\theoremseparator{:}
\newtheorem{defi}{Definition}

\theoremheaderfont{\sc}
\theorembodyfont{\upshape}
\theoremstyle{nonumberplain}
\theoremseparator{.}
\newtheorem{dem}{Proof}
    \theoremclass{dem}
    \theoremstyle{nonumberbreak}

\theoremclass{dem}
    \theoremstyle{nonumberplain}

\theoremclass{dem}
    \theoremstyle{nonumberplain}
    \theoremseparator{ }

\theoremstyle{plain}
\theoremheaderfont{\normalfont\bfseries}
\theorembodyfont{\slshape}
\theoremseparator{:}
\newtheorem{obs}[teo]{Remark}

    \theoremclass{obs}
    \theoremstyle{break}

\title{The Fundamental Theorem of Tropical Differential Algebraic Geometry}
\author{Fuensanta Aroca\textsuperscript{a}, Cristhian Garay\textsuperscript{b}, Zeinab Toghani\textsuperscript{c}\\ 
\textsuperscript{a} fuen@im.unam.mx, \\
Universidad Nacional Aut\'onoma de M\'exico.\\ 
\textsuperscript{b} cristhian.garay@imj-prg.fr,\\
Institut de  Math\'ematiques  de Jussieu--Paris  Rive  Gauche.\\
\textsuperscript{c} zeinab.toghani@im.unam.mx,\\
Universidad Nacional Aut\'onoma de M\'exico.\\}

\begin{document}

\maketitle
\let\thefootnote\relax\footnotetext{This research was supported by PAPIIT IN104713 and IN108216, ECOS NORD M14M03 and CONACYT. Part of this work was done while the first and the third authors were at the CNRS unit\'e UMR7373 12M in Marseille, the first with a CNRS contract and the third with the support of LAISLA.

\emph{2010 Mathematics Subject Classification} 13N99, 13P15, 14T99.}
\let\thefootnote\svthefootnote
\begin{abstract}
Let $I$ be an ideal of the ring of Laurent polynomials $K[x_1^{\pm1},\ldots,x_n^{\pm1}]$ with coefficients in a real-valued field $(K,v)$. The fundamental theorem of tropical algebraic geometry states the equality $\text{trop}(V(I))=V(\text{trop}(I))$ between the tropicalization $\text{trop}(V(I))$ of the closed subscheme $V(I)\subset (K^*)^n$ and the tropical variety $V(\text{trop}(I))$ associated to the tropicalization of the ideal $\text{trop}(I)$.

In this work we prove an analogous result for a differential ideal $G$ of the ring of differential polynomials $K[[t]]\{x_1,\ldots,x_n\}$, where $K$ is an uncountable algebraically closed field of characteristic zero. We define the tropicalization $\text{trop}(\text{Sol}(G))$ of the set of solutions $\text{Sol}(G)\subset K[[t]]^n$ of $G$, and the  set of solutions $\text{Sol}(\text{trop}(G))\subset\mathcal{P}(\dob{Z}_{\geq0})^n$ associated to the tropicalization of the ideal $\text{trop}(G)$. These two sets are linked by a tropicalization morphism $\text{trop}:\text{Sol}(G)\longrightarrow \text{Sol}(\text{trop}(G))$.

We show the equality $\text{trop}(\text{Sol}(G))=\text{Sol}(\text{trop}(G))$, answering a question raised by D. Grigoriev earlier this year.
\end{abstract}

\textbf{Key words:} Differential Algebra, Tropical Geometry, Arc spaces, Power series solutions of ODE.
\section{Introduction}

The first proof of the Fundamental Theorem of Tropical Algebraic Geometry appeared in 2003 in a preprint by Einsleider, Kapranov and Lind. The proof is done for hypersurfaces and was published in 2006 \cite{EinsiedlerKapranovLind:2006}. Later it was established (in full
generality)  in \cite{SpeyerSturmfels:2004}.  Extensions to arbitrary codimension ideals and arbitrary valuations have been done subsequently (see for example \cite{ArocaIlardiLopez:2010,JensenMarkwig:2008,Aroca:2010}).

The tropical variety of a hypersurface is dual to a subdivision of the Newton polyhedron of its defining function. The Newton Polygon was introduced by Puiseux in 1850 \cite{Puiseux:1850} for plane algebraic curves and extended to differential polynomials by Fine in 1889 \cite{Fine:1889}. Both the extensions of the polygon and the polyhedron have served to prove existence theorems and to construct algorithms that compute solutions (see for example \cite{GrigorievSinger:1991,Cano:1993,ArocaCano:2001,ArocaCanoJung:2003}).

In \cite{Grigoriev:2015}, Grigoriev introduces the notion of {\em tropical linear differential equations} in $n$ variables and designs a polynomial complexity algorithm for solving systems of tropical linear differential equations in one variable. In the same preprint, Grigoriev suggests several lines for further research. One of his questions is whether a theorem such as the fundamental theorem of tropical algebraic geometry holds in this context.

More precisely, Grigoriev notes that, for a differential ideal $G$ in $n$ independent variables, we have the inclusion $\text{trop} (Sol (G))\subset Sol (\text{trop} (G))$ and asks:

`` Is it true that for any differential ideal $G$  and a family $S_1,\ldots , S_n\subset\dob{Z}_{\geq0}$ being a solution of the tropical differential equation $\text{trop}(g)$ for any $g\in G$, there exists a power series solution of $G$ whose tropicalization equals $S_1,\ldots , S_n$? ''

Here, we give a positive answer to this question when $G$ is a differential ideal of differential polynomials over the ring of formal power series $K[[t]]$,  $K$ being an uncountable algebraically closed field of characteristic zero. Our proof uses techniques developed in the theory of arc spaces (see \cite{Nash:1995}).

 In Section \ref{DiferentialAlgebraicGeometry}, the basic definitions of  differential algebraic geometry are recalled. In Sections \ref{TropSer}, \ref{TropDif} and \ref{Tropicalization}, we explain the tropicalization morphisms. Arc spaces and their connection with sets of solutions of differential ideals are discussed in Section \ref{ArcSpaces}. The main result is proved in the last two sections.

The authors would like to thank Jos\'e Manuel Aroca, Ana Bravo and Guillaume Rond for fruitful conversations during the preparation of this work. The second author would also like to thank Andrew Stout, and the third author wishes to thank Pedro Luis de Angel.

\section{Differential Algebraic Geometry}\label{DiferentialAlgebraicGeometry}

We will begin by recalling some basic definitions of  differential algebraic geometry. The  reference for this section is the book \cite{Ritt:1950} by J. F. Ritt.

Let $R$ be a commutative ring with unity. A \enf{derivation} on $R$ is a map $d:R\longrightarrow R$ that satisfies   $d(a+b)=d(a)+d(b)$ and $d(ab)=d(a)b+ad(b)$,$\forall a,b\in R$. The pair $(R,d)$ is called a \enf{differential ring}. An ideal $I\subset R$ is said to be a \enf{differential ideal} when $d(I)\subset I$.

Let $(R,d)$ be a  differential ring and let $R\{x_1,\ldots,x_n\}$ be the set of polynomials with coefficients in $R$ in the variables $\{x_{ij}\::\:i=1,\ldots,n,\:j\geq0\}$. The derivation $d$ on $R$ can be extended to a derivation $d$ of $R\{x_1,\ldots,x_n\}$ by setting $d(x_{ij})=x_{i(j+1)}$ for  $i=1,\dots,n$ and $j\geq0$. The pair $(R\{x_1,\ldots,x_n\},d)$ is a differential ring called \enf{the ring of differential polynomials in $n$ variables with coefficients in $R$}.

A differential polynomial $P\in R\{x_1,\ldots,x_n\}$ induces a mapping from $R^n$ to $R$ given by 
\begin{equation}\label{SingificadoDeSubstituir}
	\begin{array}{cccc}
		P:
			& R^n
				&\longrightarrow 
					&R\\
					\\
			&(\varphi_1,\ldots,\varphi_n)
				&\mapsto
					& \displaystyle
					P|_{x_{ij}=d^j\varphi_i},
	\end{array}
\end{equation}
where $P|_{x_{ij}=d^j\varphi_i}$ is the element of $R$ obtained by substituting $x_{ij}\mapsto d^j\varphi_i$ in the differential polynomial $P$.

The equality
\begin{equation}\label{ParaFormula1}
	d^k\left( P \left(\varphi\right)\right) = \left(d^k P \right)\left(\varphi \right)
\end{equation}
holds for any $P\in R\{x_1,\ldots,x_n\}$ and $\varphi\in R^n$.

 A \enf{zero} or a \enf{solution} of $P\in R\{x_1,\ldots,x_n\}$ is an $n$-tuple $\varphi\in R^n$ such that $P(\varphi)=0$. An $n$-tuple $\varphi\in R^n$ is a \enf{solution} of $\Sigma\subset R\{x_1,\ldots,x_n\}$ when it is a solution of every differential polynomial in $\Sigma$; that is,
\[
	Sol (\Sigma):=\{\varphi\in R^n : P(\varphi)=0,\forall P\in \Sigma\}.
\]

The following result can be found in Ritt's book \cite[p. 21]{Ritt:1950}.

\begin{p}\label{Rit}\label{numero finito de ecuaciones}
The solution of any infinite system of differential polynomials $\Sigma\subset F\{x_1,\ldots,x_n\}$, where $F$ is a differential field of characteristic zero, is the solution of some finite subset  
of the system.
\end{p}

A \enf{differential monomial} in $n$ independent variables of order less than or equal to $r$ is an expression of the form
\begin{equation}\label{dm}
E_M := \prod_{\substack{1\leq i\leq n\\ 0\leq j\leq r}}x_{ij}^{M_{ij}},
\end{equation}
where $M=(M_{ij})_{\substack{1\leq i\leq n\\ 0\leq j\leq r}}$ is a matrix in $\mathcal{M}_{n\times(r+1)}(\dob{Z}_{\geq0}).$

 With this notation, a differential polynomial $P\in R\{x_1,\ldots,x_n\}$ is an expression of the form
\begin{equation}\label{pdif2}
	P=\sum_{M\in\Lambda\subset \mathcal{M}_{n\times(r	+1)}(\dob{Z}_{\geq0})}\psi_ME_M,
\end{equation}
with $r\in \dob{Z}_{\geq0}$, $\psi_M\in R$ and $\Lambda$ finite.

The mapping induced by the monomial $E_M$ is given by 
\[
	\begin{array}{cccc}
		E_M:
			& R^n
				&\longrightarrow 
					&R\\
					\\
			&(\varphi_1,\ldots,\varphi_n)
				&\mapsto
					& \displaystyle\prod_{\substack{1\leq i\leq n\\ 0\leq j\leq r}}{(d^j\varphi_i)}^{M_{ij}}.
	\end{array}
\]
and the map (\ref{SingificadoDeSubstituir}) induced by the differential polynomial $P$ in \eqref{pdif2} is 
\begin{equation}\label{MapeoPolinolioDif}
	\begin{array}{cccc}
		P:
			& R^n
				&\longrightarrow 
					&R\\
					\\
			&\varphi=(\varphi_1,\ldots,\varphi_n)
				&\mapsto
					& \displaystyle\sum_{M\in\Lambda}\psi_ME_M(\varphi).
	\end{array}
\end{equation}

\section{The differential ring of formal power series and tropicalization}\label{TropSer}

In what follows, we will work with the differential valued ring $R=K[[t]]$ where $K$ is an uncountable algebraically closed field of characteristic zero. We will denote $F=\text{Frac}(R)$. 

The elements of $R$ are expressions of the form
\begin{equation}\label{ExpresionElementoR}
	\varphi = \sum_{j\in\dob{Z}_{\geq0}} a_j t^j
\end{equation}
with $a_j\in K$ for $j\in\dob{Z}_{\geq0}$.

The \enf{support}  of $\varphi$ is the set 
\[
\text{Supp}(\varphi) :=\{i\in\dob{Z}_{\geq0}\::\:a_i\neq0\},
\]
the valuation on $R$ is given by
 \[
\text{val} (\varphi ) = \min \text{Supp}(\varphi)
\]
and the derivative of $\varphi$ is the element of $R$ 
\[
d\varphi = \sum_{j\in\dob{Z}_{\geq0}} j a_j t^{j-1}.
\]

The bijection between $K^{\dob{Z}_{\geq 0}}$ and $R$
\[
\begin{array}{cccc}
\Psi :
	&K^{\dob{Z}_{\geq 0}}
		& \longrightarrow
			&  R\\
	& \underline{a}=(a_{j})_{\substack{j\geq 0}} 
		&\mapsto
			&\displaystyle  \sum_{j\geq 0}\frac{1}{j!}a_{j}t^j
\end{array}
\]
allows us to identify points of $R$ with points of $K^{\dob{Z}_{\geq 0}}$.

 Moreover, the mapping $\Psi$ has the following property:
\begin{equation}\label{DerivadaSerie}
 d^s \Psi (\underline{a}) = \sum_{j\geq0}\frac{a_{s+j}}{j!}t^j
 \end{equation}
 which implies
 \[
 d^s \Psi (\underline{a})|_{t=0} = a_s
 \]
 and then
 \begin{equation}\label{Recuperaraevaluandoencero}
 \underline{a}=(\left. d^j \Psi ( \underline{a}) \right|_{t=0})_{ j\geq 0}.
\end{equation}

 The mapping that sends each series in $R$ to its support set (a subset of $\dob{Z}_{\geq 0}$) will be called the \enf{tropicalization} map
\[
	\begin{array}{cccc}
		\text{trop}:
			& R
				&\longrightarrow
					& \mathcal{P}(\dob{Z}_{\geq 0})\\
			& \varphi
				&\mapsto
					&\text{Supp}(\varphi)
	\end{array}
\]
where $\mathcal{P}(\dob{Z}_{\geq0})$ denotes the power set of $\dob{Z}_{\geq0}$.

For fixed $n$, the mapping  from $R^n$ to the $n$-fold product of $\mathcal{P}(\dob{Z}_{\geq 0})$ will also be denoted by $\text{trop}$:
\[
	\begin{array}{cccc}
		\text{trop}:
			& R^n
				&\longrightarrow
					& {\mathcal{P}(\dob{Z}_{\geq 0})}^n\\
			& \varphi=(\varphi_1,\ldots ,\varphi_n)
				&\mapsto
					& \text{trop}(\varphi)=(\text{Supp}(\varphi_1),\ldots ,\text{Supp}(\varphi_n)).
	\end{array}
\]

Given a subset $T$ of $R^n$, the \enf{tropicalization} $T$ is its image under the map $\text{trop}$:
\[
\text{trop} (T):=\left\{ \text{trop}(\varphi) \::\: \varphi\in T\right\}\subset {\mathcal{P}(\dob{Z}_{\geq 0})}^n.
\]
\begin{ej}

Set $T:= \{ ( a+5t+bt^2, 2+ at -8t^2 + ct^3) : a,b,c\in K\}\subset {K[[t]]}^2$ we have
\[
\begin{array}{ll}
	\text{trop}(T) = 
		&\{ (\{ 1\}, \{ 0,2\}), (\{0, 1\}, \{ 0,1,2\}), (\{ 1,2\}, \{ 0,2\}),(\{ 1\}, \{ 0,2,3\}),\\
		&(\{0, 1,2\}, \{ 0,1,2\}),(\{0, 1\}, \{ 0,1,2,3\}), (\{ 1,2\}, \{ 0,2,3\}),(\{ 0,1,2\}, \{ 0,1,2,3\} ) \}.
\end{array}
\]
\end{ej}

Since $K$ is of characteristic zero, for every $\varphi\in R$, we have
\[
\text{trop} \left( d^j\varphi\right) = \left\{ i-j : i\in \text{trop} (\varphi)\cap \dob{Z}_{\geq j}\right\}
\]
then
\[
\text{val}\left( d^j\varphi\right) = \min \left(\text{trop} (\varphi)\cap \dob{Z}_{\geq j}\right) -j.
\]

The above equality justifies the following definition:
\begin{defi}				
A subset $S\subseteq\dob{Z}_{\geq0}$ induces a mapping   $\text{Val}_S:\dob{Z}_{\geq0}\longrightarrow\dob{Z}_{\geq0}\cup\{\infty\}$ given by 
\begin{equation}\label{tropsol2}
\text{Val}_S(j):= \begin{cases}
s-j,&\text{ with }s=\text{min}\{\alpha\in S\::\:\alpha\geq j\},\\
\infty,&\text{ when } S\cap \dob{Z}_{\geq j}=\emptyset.
\end{cases}
\end{equation}
\end{defi}

\begin{ej} Consider the set $S:=\{1,3,4\}$. We have
\begin{enumerate}
\item  $\text{Val}_{S}(2)=\min \{s\in S \mid s\ge 2\}-2=3-2=1$
\item $\text{Val}_{S}(5)=\infty$.
\end{enumerate}
\end{ej}

\section{Tropical differential polynomials}\label{TropDif} 

We will denote by $\dob{T}$ the (idempotent) semiring $\dob{T}=(\dob{Z}_{\geq0}\cup\{\infty\},\oplus,\odot)$, with $a\oplus b=\text{min}\{a,b\}$ and $a\odot b=a+b$. 
\begin{defi}
A \enf{tropical differential monomial} in the variables $x_1,\ldots,x_n$ of order less than or equal to $r$ is an expression of the form
\begin{equation}\label{tdm2}
\varepsilon_M:=x^{\odot M}=\bigodot_{\substack{1\leq i\leq n\\ 0\leq j\leq r}}x_{i,j}^{\odot M_{ij}},
\end{equation}
where $M=(M_{ij})_{\substack{1\leq i\leq n\\ 0\leq j\leq r}}$ is a matrix in $\mathcal{M}_{n\times(r+1)}(\dob{Z}_{\geq0})$.
\end{defi}

\begin{defi}
A \enf{tropical differential polynomial}  in the variables $x_1,\ldots,x_n$ of order less than or equal to $r$ is an expression of the form 
\begin{equation}\label{tdp}
\phi=\phi(x_1,\ldots,x_n)=\underset{M\in\Lambda\subset \mathcal{M}_{n\times(r+1)}(\dob{Z}_{\geq0})}{\bigoplus}a_{M}\odot\varepsilon_M,
\end{equation}
where $a_{M}\in\dob{T}$ and $\Lambda$ is a finite set.
\end{defi}

The set of tropical differential polynomials will be  denoted by $\dob{T}\{x_1,\ldots,x_n\}$. 

A tropical differential monomial $\varepsilon_M$   induces a mapping from $\mathcal{P}(\dob{Z}_{\geq0})^n$ to $\dob{Z}_{\geq0}\cup\{\infty\}$ given by 
\begin{equation*}
\varepsilon_M(S_1,\ldots,S_n):=  \bigodot_{\substack{1\leq i\leq n\\ 0\leq j\leq r}}{\text{Val}_{S_i}(j)}^{\odot M_{ij}} = \sum_{\substack{1\leq i\leq n\\ 0\leq j\leq r}}M_{ij}\cdot \text{Val}_{S_i}(j),
\end{equation*}
where  $\text{Val}_{S_i}(j)$ is defined as in \eqref{tropsol2}.

\begin{obs}\label{ParaElResultado}
Note that  $\varepsilon_M(S_1,\ldots,S_n)=0$ if and only if $j\in S_i$ for all $i, j$ with $M_{ij}\neq0$.
\end{obs}

 A tropical differential polynomial $\phi$ as in \eqref{tdp} induces a mapping from $\mathcal{P}(\dob{Z}_{\geq0})^n$ to $\dob{Z}_{\geq0}\cup\{\infty\}$ given by 
 \[
 \phi(S)=\underset{M\in\Lambda}{\bigoplus}a_{M}\odot\varepsilon_M(S)=\underset{M\in\Lambda}{\text{min}}\{a_{M}+\varepsilon_M(S)\}.
\]

\begin{defi}
An $n$-tuple $S=(S_1,\ldots,S_n)\in \mathcal{P}(\dob{Z}_{\geq0})^n$ is said to be a \enf{solution} of  the tropical differential polynomial $\phi$ in \eqref{tdp}   if either
\begin{enumerate}
\item There exists $M_1,M_2\in\Lambda$, $M_1\neq M_2$, such that $\phi(S)=a_{M_1}\odot\varepsilon_{M_1}(S)=a_{M_2}\odot\varepsilon_{M_2}(S)$, or
\item $\phi(S)=\infty$.
\end{enumerate}
\end{defi}
Let $H\subset\dob{T}\{x_1,\ldots,x_n\}$ be a family of tropical differential polynomials. An $n$-tuple $S\in \mathcal{P}(\dob{Z}_{\geq0})^n$ is a \enf{solution} of $H$ when it is a solution of every tropical polynomial in $H$; that is,
\[
Sol (H) := \left\{ S\in {(\mathcal{P}(\dob{Z}_{\geq0}))}^n : S \text{ is a solution of } \phi \text{ for every }\phi\in H\right\}.
\]

\begin{ej}
Consider the tropical differential polynomial
\[
\phi(x):= 1\odot x' \oplus  2\odot x^{(3)}  \oplus 3.
\]
Since $\phi(S)\neq\infty$ for every $S\subset\mathcal{P}(\dob{Z}_{\geq0})$, the set $S$ is a solution of $\phi$ if one of the following holds
\begin{enumerate}
\item $1+\text{Val}_S(1)=3 \le 2+\text{Val}_S(3)$,
\item $1+\text{Val}_S(1)= 2+\text{Val}_S(3)\le 3$,
\item $2+\text{Val}_S(3)=3 \le 1+\text{Val}_S(1)$.
\end{enumerate}
The first condition never holds.
The second condition holds when $ S=B\cup\{2,3\}\cup C$ and $B\subset\{0\}$, $\min C \ge 4$.
The third condition holds when $ S=\{4\}\cup C\cup B$ with $\min C \geq 5$ and  $B\subset \{0\}$.
\[
 Sol(P)=\{B\cup\{2,3\}\cup C ;B\subset\{0\},\:\min C \ge 4 \}\cup \{B\cup\{4\}\cup C ;  \min C \ge 5  ,  B\subset \{0\}\}.
 \]
\end{ej}

\section{Tropicalization of differential polynomials} \label{Tropicalization}

Let $P$ be a differential polynomial as in Equation \eqref{pdif2}. The \enf{tropicalization} of $P$ is the tropical differential  polynomial 
\begin{equation}\label{tdpK}
\text{trop}(P): =\underset{M\in\Lambda}{\bigoplus}\text{val}(\psi_{M})\odot\varepsilon_M.
\end{equation}

\begin{obs}\label{ParaElResultado2}
Let $P$ be a differential polynomial in $R\{ x_1,\ldots ,x_n\}$.  We have that $\mathrm{trop}(tP )(S)\geq 1$ for any $S\in {\cal{P}} (\dob{Z}_{\geq 0})^n$.
\end{obs}

\begin{defi}
Let $G\subset R\{x_1,\ldots,x_n\}$ be a differential ideal. Its \enf{tropicalization} $\text{trop}(G) $ is the set of tropical differential  polynomials $\{ \text{trop}(P) : P\in G\}$.
\end{defi}

\begin{p}\label{tropofsol}
Let $G$ be a differential ideal in the ring of differential polynomials $R\{x_1,\ldots,x_n\}$. If $\varphi\in\text{Sol}(G)$, then $\text{trop}(\varphi)\in Sol (\text{trop}(G))$.
\end{p}
\begin{dem}
Given a differential monomial $E_M$ and $\varphi\in R^n$, we have that $$\text{val}(E_M(\varphi))=\varepsilon_M(\text{trop}(\varphi)).$$
It follows that if $\varphi\in R^n$ is a solution to the differential polynomial $P=\sum_{M\in\Lambda}\psi_ME_M$, then $\text{trop}(\varphi)\in(\mathcal{P}(\dob{Z}_{\geq0}))^n$ is a solution to $\text{trop}(P)$. So, if $\varphi\in R^n$ is a solution to every differential polynomial $P$ in a differential ideal $G$, then $\text{trop}(\varphi)$ is a solution to every tropical differential polynomial $\text{trop}(P)\in\text{trop}(G)$.
\fin
\end{dem}

We can now clearly state the question  posed by Grigoriev in \cite{Grigoriev:2015}. The latter result allows us to define a mapping $\text{trop}: \text{Sol} (G)\longrightarrow \text{Sol}(\text{trop}(G))$ for any differential ideal $G\subset R\{x_1,\ldots,x_n\}$. The question is whether or not this map is surjective.
\begin{ej}
 Let $P\in R\{x\}$ be the differential polynomial
 \[
 	P:= x''-t.
 \]
 The set of solutions of $P$ is the same as the set of solutions of the differential ideal generated by $P$
 \[
 \text{Sol}(P) = \{ c_1 +c_2 t +\frac{1}{6}t^3 : c_1,c_2\in K\}.
 \]
 The tropicalization of the set of solutions of $P$ is
 \[
 \text{trop}(\text{Sol} (P)) = \left\{ \{ 0,1,3\}, \{ 0,3\},\{ 1,3\}, \{ 3\}
 	\right\}.
 \]
 Now, the tropicalization of $P$ induces the mapping
 \[
 	\begin{array}{cccc}
 		\text{trop}(P):
 			&\mathcal{P}(\dob{Z}_{\geq0})
 				&\longrightarrow
 					&\dob{Z}_{\geq0}\\
 			& S	&\mapsto
 					& \text{min}\{\text{Val}_S(2),1\}.
 	
 	\end{array}
 \]
 Since $\text{trop}(P)(S)\neq\infty$ for every $S\subset\mathcal{P}(\dob{Z}_\geq0)$, the set of solutions of $\text{trop}(P)$ is
 \[
 	\text{Sol}(\text{trop}(P)) = \left\{ S\subset\mathcal{P}(\dob{Z}_\geq0)\::\:2\notin S\text{ and } 3\in S
 	\right\}.
 \]
 Differentiating $P$, we have that $d^2P= x^{(4)}$ is in the differential ideal generated by $P$. Its tropicalization induces the mapping
 \[
 	\begin{array}{cccc}
 		\text{trop}(d^2P):
 			&\mathcal{P}(\dob{Z}_{\geq0})
 				&\longrightarrow
 					&\dob{Z}_{\geq0}\\
 			& S	&\mapsto
 					& \text{Val}_S(4).
 	\end{array}
 \]
We have that $S\subset\mathcal{P}(\dob{Z}_\geq0)$ is a solution of  $\text{trop} (d^2P)$ if and only if $\text{S}\subset\{0,1,2,3\}$, i.e.
 \[
 	\text{Sol} (\text{trop} (d^2P)) = \mathcal{P}(\{0,1,2,3\}).
 \]
In this example
 \[
 \text{Sol} (\text{trop} (P)) \cap \text{Sol} (\text{trop} (d^2P))= \text{trop} (\text{Sol} (P)).
 \]
\end{ej}

\section{Arc spaces and the set of solutions of a differential ideal}\label{ArcSpaces}

The natural inclusion $K [x_{1\nobreak\hspace{.1em} 0},\ldots ,x_{n\nobreak\hspace{.1em} 0}]\subset R\{x_1,\ldots,x_n\}$ lets us recognize the arc space of the variety defined by an ideal $I \subset K [x_1,\ldots ,x_n]$ as the space of solutions of the differential ideal generated by $I$ in $R\{x_1,\ldots,x_n\}$. In this section we extend some definitions and results developed in the theory of arc spaces (see for example  \cite{Nash:1995,BruschekMourtadSchepers:2013}).

Consider the bijection 

\[
\begin{array}{cccc}
\Psi :
	&(K^{\dob{Z}_{\geq 0}})^n
		& \longrightarrow
			&  R^n\\
	& \underline{a}=(a_{ij})_{\substack{1\leq i\leq n\\ j\geq 0}} 
		&\mapsto
			&\displaystyle \left( \sum_{j\geq 0}\frac{1}{j!}a_{1j}t^j,\ldots ,\sum_{j\geq 0}\frac{1}{j!}a_{nj}t^j \right).
\end{array}
\]

\begin{lema}\label{FormulaEnCoeficientes}
Given $P\in R\{x_1,\ldots,x_n\}$ and  $\underline{a}\in (K^{\dob{Z}_{\geq 0}})^n$, we have
 \begin{equation}\label{PdePsi}
 	P(\Psi( \underline{a}) )=\sum_{k\geq 0} c_k t^k 
 \end{equation}
 with
 \[
c_k =\frac{1}{k!}\left.(d^k(P))\right|_{t=0}( \underline{a} ).
 \]
\end{lema}
\begin{dem}
 For $\underline{a}=(a_{ij})_{\substack{1\leq i\leq n\\ j\geq 0}} \in (K^{\dob{Z}_{\geq 0}})^n$, write $\Psi ( \underline{a})=(\Psi ( \underline{a})_1,\ldots,\Psi ( \underline{a})_n)$
 and
 $P(\Psi( \underline{a}) )=\sum_{k\geq 0} c_k t^k$
 for some $c_k\in K$, $k\geq 0$. Differentiating (\ref{PdePsi}) and evaluating in zero, we have 
 
 \begin{equation*}
 \begin{aligned}
 c_k=&\frac{1}{k!}[d^k(P(\Psi(\underline{a})))]_{t=0}\stackrel{(\ref{ParaFormula1})}{=} \frac{1}{k!}[(d^kP)(\Psi(\underline{a}))]_{t=0}\stackrel{(\ref{SingificadoDeSubstituir})}{=}\frac{1}{k!}
 \left[(d^kP)|_{x_{ij}=\Psi(\underline{a})_i^{(j)}}\right]_{t=0}=\\
 &\\
 &=\frac{1}{k!}
 \left[(d^kP)|_{x_{ij}=\Psi(\underline{a})_i^{(j)}|_{t=0}}\right]_{t=0}\stackrel{(\ref{Recuperaraevaluandoencero})}{=}\frac{1}{k!}\left[(d^kP)|_{x_{ij}=a_{ij}}\right]_{t=0}=\frac{1}{k!}(d^kP)|_{t=0}(\underline{a}).
 \end{aligned}
\end{equation*}

\fin
\end{dem}

Let  $G$ be a differential ideal in $R\{x_1,\ldots,x_n\}$. We can consider $G$ as an infinite system of differential polynomials in $F\{x_1,\ldots,x_n\}$, where $F=Frac(R)$ is a field of characteristic zero. By Proposition \ref{numero finito de ecuaciones}, there exist $f_1, \ldots , f_s\in G$ such that 
\[
Sol(G)=\bigcap_{\ell=1}^sSol(f_\ell).
\]
For $1\leq \ell\leq s$ and $k\in \dob{Z}_{\geq 0}$, the $\left.(d^kf_\ell)\right|_{t=0}$ are polynomials in the variables $x_{ij}$ with coefficients in $K$. Set
\begin{equation*}
 F_{\ell k}:= (d^kf_\ell)|_{t=0}\in K[ x_{ij} : 1\leq i\leq n, j\geq 0]
\end{equation*}
and
\begin{equation}\label{definicionDeAinfinito}
 A_\infty := V \left(\{ F_{\ell k}\}_{\substack{1\leq \ell\leq s\\ k \geq 0}}\right)\subset {\left( K^{\dob{Z}_{\geq 0}}\right)}^n.
\end{equation}
By Lemma \ref{FormulaEnCoeficientes}, 
\[
Sol(G)=\Psi ( A_\infty).
\]

We will now describe an extension to differential ideals of the definition of $m$-jet of arc spaces (see for example \cite{Mourtada:2011}).
 
For each   $m\geq 0$, let $N_m$ be the smallest positive integer such that
\begin{equation}\label{EleccionDeNm}
F_{\ell,k}\in K[x_{ij}, 1\leq i\leq n, 0\leq j\leq N_m],\,\forall\, 1\leq \ell\leq s,\:0\leq k\leq m
\end{equation}
and set
\begin{equation}\label{definicionDeAm}
	A_m := V \left(\{ F_{\ell k}\}_{\substack{1\leq \ell\leq s\\ 0\leq k \leq m}}\right)\subset {\left( K^{N_m+1}\right)}^n.
\end{equation}

For $m\geq m' \geq0$, denote by $\pi_{(m,m')}$ the
 natural algebraic morphism 
 \[
 \pi_{(m,m')}:K^{n(N_m+1)}\longrightarrow K^{n(N_{m'}+1)}
 \]
 then
\[
\pi_{(m,m')}(A_m)\subset A_{m'}
\]
and $A_\infty$ is the inverse limit of the system  $((A_m)_{m\in\dob{Z}_{\geq 0}}, (\pi_{(m,m')})_{m\geq m' \in\dob{Z}_{\geq 0}} )$
\[
A_\infty = \varprojlim A_m.
\]

When $f_1,\ldots ,f_s$ are elements of $K [x_{1\nobreak\hspace{.1em} 0},\ldots ,x_{n\nobreak\hspace{.1em} 0}]$ the sets $A_m$ are the $m$-jets of the space $A_\infty$. Otherwise, note that the construction depends strongly on the choice of $f_1,\ldots ,f_s$.

\section{Intersections with tori}

  Let $G\subset R\{x_1,\ldots,x_n\}$ be  a differential ideal, let $f_1,\ldots ,f_s\in G$ be such that $Sol(G)=\bigcap_{\ell=1}^sSol(f_\ell)$, and let $A_\infty$ be as in (\ref{definicionDeAinfinito}) and $A_m$ as in (\ref{definicionDeAm}).

  An $n$-tuple $S=(S_1,\ldots,S_n)\in {\mathcal{P}(\dob{Z}_{\geq 0})}^n$ is in $\text{trop}(Sol (G))$ if and only if there exists $\underline{a}\in A_\infty$ with $\text{trop} (\Psi (\underline{a} ))=S$, i.e., if $S_i=\{ j \::\: a_{ij}\neq 0\}$ for $i=1,\ldots,n$.
  
  Set
\[
{\dob{V}_S}^*:= \left\{ \left( x_{ij}\right)_{1\leq i\leq n, j\geq 0}\in {\left( K^{\dob{Z}_{\geq 0}}\right)}^n\: :\: x_{ij}=0 \text{ if and only if }  j\notin S_i \right\},
\]
then $S\in\text{trop}(Sol (G))$ if and only if
\[
{A_\infty}_S := A_\infty\cap \dob{V}_S^*
\]
is not empty.

For $m\geq0$, consider the finite dimensional torus
\[
{{\dob{V}_m}_S}^*:= \left\{ \left( x_{ij}\right)_{1\leq i\leq n, 0\leq j\leq N_m}\in K^{n(N_m+1)}\: :\: x_{ij}=0 \text{ if and only if }  j\notin S_i \right\},
\]
where $N_m$ is the minimum such that (\ref{EleccionDeNm}) holds. We have ${{\dob{V}_m}_S}^*\simeq (K^*)^{L_m}$, with $L_m\leq n(N_m +1)$.
 
   Set
 \[
{A_m}_S := A_m \cap {{\dob{V}_m}_S}^*,
\]
for $m\geq m' \geq0$, the inclusions $\pi_{(m,m')}({\dob{V}_m}_S^*)\subset {\dob{V}_{m'}}_S^*$ and $\pi_{(m,m')}({A_m}_S)\subset {A_{m'}}_S$ hold, and 
 ${A_\infty}_S$ is the inverse limit of the system  $(({{A_m}_S)}_{m\in\dob{Z}_{\geq 0}}, (\pi_{(m,m')})_{m\geq m' \in\dob{Z}_{\geq 0}})$,
 \[
{A_\infty}_S = \varprojlim {A_m}_S.
\]
Set 
\[
{B_m}_S := \displaystyle\bigcap_{i=m}^\infty \pi_{(i,m)}({A_i}_S)
\]
then
\[
{A_\infty}_S = \varprojlim {B_m}_S
\]
and the projections
\[
\pi_{(m,m')}: {B_m}_S\longrightarrow {B_{m'}}_S
\]
are surjective. Then (see for example \cite[Proposition 5, p.198]{Bourbak:2004}),  the set $\varprojlim {B_m}_S$ is nonempty if and only if ${B_0}_S$ is nonempty. In other words, we have the following remark.

\begin{obs}\label{AinfinitoVacioSiInteresccionVacia}
 The set ${A_\infty}_S$ is nonempty if and only if 
 the intersection $\displaystyle\bigcap_{i=0}^\infty \pi_{(i,0)}({A_i}_S)$ is nonempty.
\end{obs}

By Chevalley's Theorem (see for example \cite[p.51]{Mumford:1999}), each $\pi_{(m,0)}({A_m}_S)$ is a constructible set. A constructible set is, by definition, a finite union of locally closed sets. A set is locally closed when it is an open set of its closure. The constructible sets form a Boolean algebra.

We recall the following statement about nested sequences of constructible sets: 
\begin{p}\label{EGA}
  Let $ K $ be an uncountable algebraically closed field of characteristic zero. Let  $\{E_\alpha\}_{\alpha=1}^{\infty}$ be an increasing family of constructible sets in $K^n$ with $K^n=\bigcup_{\alpha=1}^{\infty} E_\alpha$, then there exists $\alpha$ such that $K^n=E_\alpha$. 
\end{p}

We are now ready to prove the result that will allow us, in the next section, to work in the noetherian ring $K[x_{ij}, 1\leq i\leq n, 0\leq j\leq N_m]$ instead of the non-noetherian $K[x_{ij}, 1\leq i\leq n, 0\leq j]$.

\begin{p}\label{SiLaInterseccionEsNOVacia}
The set ${A_\infty}_S$ is nonempty if and only if $\displaystyle {A_m}_S$ is nonempty for all $m\in \dob{Z}_{\geq 0}$.
\end{p}

\begin{dem}

  Since the constructible sets form a Boolean algebra, the nested sequence of constructible sets inside $(K^*)^{L_0}\simeq \dob{V}_{0S}^*$
\begin{equation}\label{SucesionAnidada}
\cdots\:\subset\:\pi_{(2,0)}({A_2}_S)\:\subset\:\pi_{(1,0)}({A_1}_S)\:\subset\: {A_0}_S\:\subset\: {\left( K^*\right)}^{L_0}
\end{equation}
induces an increasing family of constructible sets
\begin{equation}
\emptyset \enspace \subset \enspace {\left( K^*\right)}^{L_0}\setminus {A_0}_S \enspace \subset \enspace {\left( K^*\right)}^{L_0}\setminus \pi_{(1,0)} \left( {A_1}_S\right)\enspace \subset \enspace {\left( K^*\right)}^{L_0}\setminus \pi_{(2,0)} \left( {A_2}_S\right)\enspace \subset \enspace \cdots
\end{equation}

The set $\displaystyle\bigcap_{i=0}^\infty \pi_{(i,0)}({A_i}_S)$ is empty if and only if $\displaystyle {\left( K^*\right)}^{L_0}\setminus\bigcap_{i=0}^\infty \pi_{(i,0)}({A_i}_S)$ is $\displaystyle {\left( K^*\right)}^{L_0}$. That is, if and only if
\[
{\left( K^*\right)}^{L_0} = \bigcup_{i=0}^\infty {\left( K^*\right)}^{L_0}\setminus \pi_{(i,0)} \left( {A_i}_S\right).
\]
Then, by Proposition \ref{EGA}, there exists $m$ such that ${\left( K^*\right)}^{L_0}\setminus \pi_{(m,0)} \left( {A_m}_S\right)={\left( K^*\right)}^{L_0}$. That is, there exists $m$ such that ${A_m}_S$ is empty.

The result follows from Remark \ref{AinfinitoVacioSiInteresccionVacia}.
\fin
\end{dem}

\section{The fundamental theorem of differential tropical geometry}

\begin{teo}
Let $G$ be a differential ideal in $K[[t]]\{ x_1,\ldots x_n\}$, where $K$ is an uncountable algebraically closed field of characteristic zero. The equality
\[
\mathrm{Sol} (\mathrm{trop}(G))= \mathrm{trop}(\mathrm{Sol} (G))
\]
holds.
\end{teo}
\begin{dem}
The inclusion $  \text{trop}(\text{Sol}  (G))\subset \text{Sol}  (\text{trop}(G)) $ is just Proposition \ref{tropofsol}. Here we will prove  
\[
\text{Sol} (\text{trop}(G))\subset \text{trop}(\text{Sol}  (G)).
\]

Let $S=(S_1,\ldots,S_n)\in {\mathcal{P}(\dob{Z}_{\geq 0})}^n$ be such that there is no solution of $G$ whose tropicalization is $S$. We will show that $S$ cannot be a solution of the tropicalization of $G$.

Suppose that $\text{Sol} (G)=\bigcap_{\ell=1}^s\text{Sol} (f_\ell)$, for some $f_1,\ldots f_s\in G$. For $1\leq \ell\leq s$ and $k\in \dob{Z}_{\geq 0}$, we write $ F_{\ell k}:= (d^kf_\ell)|_{t=0}$.

As we have seen above, $ S\notin \text{trop}(\text{Sol}  (G))$ implies that $A_{\infty S}$ is empty. Then, by Proposition \ref{SiLaInterseccionEsNOVacia} there exists   $m\in\dob{N}$ such that ${A_m}_S$ is empty.

Take $m\in\dob{N}$ such  that ${A_m}_S$ is empty. Set 
$\overline{F_{\ell k}}$ to be the image of $F_{\ell k}$ in the ring $$K [ x_{ij} : 1\leq i\leq n,0\leq j\leq N_{m}]/\left<x_{ij}\::\:j\notin S_i\right>.$$

Since ${A_m}_S$ is empty we have 
\[
V\left(\overline{F_{\ell,k}}\::\:1\leq \ell\leq s,\:0\leq k\leq m\right)\subset V\left(\prod_{\{0\leq i\leq n, j\in S_i\::\:j\leq N_m\}}x_{i,j}\right)
\]
 so by the Nullstellensatz, there exists $\alpha\geq 1$ such that 
\[ 
 E_M={\prod_{\{0\leq i\leq n, j\in S_i\::\:j\leq N_m\}}x_{i,j}}^\alpha
\in\left<\overline{F_{\ell,k}
 }\::\:1\leq \ell\leq s,\:0\leq k\leq m\right>.
\]
Here $E_M$ is the differential monomial induced by the matrix $M\in\mathcal{M}_{n\times(N_m+1)}(\dob{Z}_{\geq0})$  with entries $M_{ij}=0$ for $j\notin S_i$ and $M_{ij}=\alpha$ for $j\in S_i$. 

It follows that there exists 
 $\{G_{\ell,k}\::\:1\leq \ell\leq s,\:0\leq k\leq m\}\subset K[x_{ij}, 1\leq i\leq n, j\in S_i, j\leq N_m]$
such that \begin{equation}
 \sum_{\substack{1\leq \ell\leq s\\ 0\leq k\leq m}}G_{\ell,k}\overline {F_{\ell,k}}=E_M.
 \end{equation}
 then
  \begin{equation}
 \sum_{\substack{1\leq \ell\leq s\\ 0\leq k\leq m}}G_{\ell,k}F_{\ell,k}=E_M+h
 \end{equation}
for some  $h\in \left<x_{ij}\::\:j\notin S_i,\:j\leq N_m\right>\subset K [ x_{ij} : 1\leq i\leq n,0\leq j\leq N_{m}]$.

 Now, by definition of $F_{\ell k}$,   there exists $\lambda$ in $K[[t]]\{x_0,\ldots ,x_n\}$ so that 
 \begin{equation}
 g:=\sum_{\substack{1\leq \ell\leq s\\ 0\leq k\leq m}}G_{\ell,k} d^kf_\ell=E_M+h+t\lambda .
 \end{equation}
 Since $G$ is a differential ideal and $f_1,\ldots f_s\in G$, the differential polynomial $g$ is in $G$.
 
 Now
 \begin{enumerate}
 	\item
 	By Remark \ref{ParaElResultado}
 	\begin{enumerate}
 		\item
 		$\varepsilon_M(S) =0$
 		\item
 		If $h\neq0$, then $\text{trop} (h)(S)\geq 1$.
 	\end{enumerate}
 	\item
 	By Remark \ref{ParaElResultado2}
 	\begin{enumerate}
 		\item
 		If $t\lambda\neq0$, then $\text{trop}(t\lambda)(S)\geq 1$.
 	\end{enumerate}
\end{enumerate}
 Then $(\text{trop}(g))(S)=0$ and the minimum is attained only at the momomial $\varepsilon_M$, and then, $S$ is not a solution of $\text{trop}(g)$. So $S$ is not a solution of the tropicalization of $G$ which is what we wanted to prove. \fin
 \end{dem}

\bibliographystyle{amsalpha}
\addcontentsline{toc}{chapter}{Bibliography}
\bibliography{bibliografia_zeinab}

\end{document}